\documentclass{elsart3-1}

\usepackage{amssymb}
\usepackage{amsmath}
\usepackage[english,francais]{babel}

\newtheorem{theoreme}{Th\'eor\`eme}[section]

\newtheorem{proposition}[theoreme]{Proposition}
\newtheorem{corollaire}[theoreme]{Corollaire}
\newtheorem{definition}[theoreme]{D\'efinition\rm}

{\nopagebreak[0]\hfill$\Box$ \medskip}%
\renewcommand{\theequation}{\arabic{equation}}
\def\N{\hbox{\bf N}}

\def\F{\hbox{\bf F}}

\def\Z{\hbox{\bf Z}}

\def\ds{\displaystyle}
\def\iy{\infty}
\def\mk{\smallskip}

\def\n{\noindent}

\def\g{\gamma}

\def\D{\Delta}

\def\f{\varphi}

\def\cf{\mathcal{F}}

\def\ni{\noindent}

\def\Lr{\Longrightarrow}

\def\numero{n$^{\text{o}}$}

\def\leq{\leqslant}
\def\geq{\geqslant}
\def\cleq{\preccurlyeq}
\def\clet{\prec}

\def\cget{\succ}

\usepackage[T1]{fontenc}

\def\og{\leavevmode\raise.3ex\hbox{$\scriptscriptstyle\langle\!\langle$~}}
\def\fg{\leavevmode\raise.3ex\hbox{~$\!\scriptscriptstyle\,\rangle\!\rangle$}}

\journal{the Acad\'emie des sciences}
\begin{document}

\centerline{}
\begin{frontmatter}




%
\selectlanguage{francais}
\title{Formes modulaires modulo $2$ : l'ordre de nilpotence des 
op\'erateurs de Hecke}


\author[authorlabel1]{Jean-Louis NICOLAS},
\ead{jlnicola@in2p3.fr,\;\;http://math.univ-lyon1.fr/$\sim$nicolas/}
\author[authorlabel2]{Jean-Pierre SERRE}
\ead{jpserre691@gmail.com}

\address[authorlabel1]{CNRS, Universit\'e de Lyon, Institut Camille
  Jordan, Math\'ematiques,  F-69622 Villeurbanne Cedex, France.}
\address[authorlabel2]{Coll\`ege de France, 3 rue d'Ulm, F-75231 Paris Cedex 05, France.}


\medskip

\selectlanguage{english}
\begin{abstract}




\noindent
{\bf The nilpotence order of the mod 2 Hecke operators.}

\vskip 0.5\baselineskip

Let $\D= \sum_{m=0}^\iy
q^{(2m+1)^2} \in \F_2[[q]]$ be the reduction mod 2 of the $\D$ series.
A modular form $f$ modulo $2$ of level 1 is a polynomial in $\D$. If $p$ is an odd prime, then the Hecke
operator $T_p$ transforms $f$ in a modular form $T_p(f)$ which is a
polynomial in $\D$ whose degree is smaller than the degree of $f$, so
that $T_p$ is nilpotent.

The order of nilpotence of $f$ is defined as the smallest integer
$g=g(f)$ such that, for every family of $g$ odd primes $p_1,p_2,\ldots,p_g$, the relation
$T_{p_1}T_{p_2}\ldots T_{p_g}(f)=0$ holds. We show how one can compute
explicitly $g(f)$; if $f$ is a polynomial of degree $d$  in $\D$, one finds that
$g(f)  < \! < d^{1/2}$.

\bigskip

\ni

\medskip

\ni
\end{abstract}
\end{frontmatter}

\selectlanguage{francais}

\section{Introduction}

Soit 
$\ds \D(q)=q\prod_{n=1}^\iy (1-q^n)^{24}=\sum_{n=1}^\iy \tau(n) q^n$ 
o\`u $\tau$ est la fonction de Ramanujan.
Soit $k$ un entier $\geq 0$. On \'ecrit
$\ds \D^k(q)=\sum_{n=k}^\iy \tau_k(n)q^n.$
Les congruences connues sur $\tau(n) \pmod{2}$ (cf. \cite{SW1}), montrent
que $\ds \D(q)\equiv \sum_{m=0}^\iy q^{(2m+1)^2} \pmod{2}$,
ce qui entra\^ine
\begin{equation}\label{cnk}
n\not\equiv k\pmod{8}\quad \Lr \quad \tau_k(n)\equiv 0 \pmod{2}.
\end{equation}
Une forme modulaire modulo $2$ de niveau 1 est un polyn\^ome $f(\D)$  \`a  coefficients
dans $\F_2$ (cf. par exemple \cite{Nic,Ser5});
nous l'identifierons  \`a  une s\'erie formelle en la variable $q$,  \`a 
coefficients dans $\F_2$. Nous ne nous int\'eresserons qu'aux formes
paraboliques (celles dont le terme constant est 0).
 \`A  partir de maintenant (sauf mention expresse du contraire), 
toutes les s\'eries
consid\'er\'ees sont  \`a  coefficients mod 2, et nous nous 
permettrons d'\'ecrire
\begin{equation}\label{Delta}
\D=\D(q)=\sum_{m=0}^\iy q^{(2m+1)^2} \  \in \F_2[[q]].
\end{equation} 

\section{Pr\'eliminaires}\label{formod}

\subsection{Les $\mathbf{\F_2}$-espaces vectoriels 
$\mathbf{\cf, \cf_1,\cf_3,\cf_5,\cf_7}$}\label{parFi}

Soit $\cf$ le sous-espace 
de $ \F_2[\D]$ engendr\'e par $\D,\D^3,\D^5,\ldots$.
Compte tenu de (\ref{cnk}), on a 
$\cf=\cf_1\oplus\cf_3\oplus\cf_5\oplus\cf_7$
o\`u, pour $i\in\{1,3,5,7\}$, $\cf_i$ a pour base $\{\D^i,\D^{i+8},
\D^{i+16},\ldots \}.$ 

Puisque $\D^{2k}(q)=\D^k(q^2)$, toute forme parabolique $f$ modulo $2$ peut
s'\'ecrire comme une somme finie 
\begin{equation}\label{fs}
f=\sum_{s\geq 0} f_s^{2^s} \quad {\rm avec} \quad f_s\in \cf.
\end{equation}

\subsection{Op\'erateurs de Hecke}
 
Soit $f(q)=\sum_{n\geq 0} c_n q^n$ une forme modulaire modulo $2$ et
soit $p$ un nombre premier $> 2$. L'op\'erateur de Hecke 
$T_p$ transforme $f$ 
en la forme 
\begin{equation}\label{ga2}
T_p|f=\sum_{n\geq 0} \g_n q^n \;\; {\rm avec } \;\;
\g (n)=\left\{
\begin{array}{ll}
c(pn)& {\rm si} \;p \; {\rm ne \; divise \; pas\; }  n\\
c(pn)+c(n/p)\;\;& {\rm si} \;p \;{\rm divise}\;  n.
\end{array}
\right.
\end{equation}
[Nous \'ecrirons parfois $T_p(f)$  \`a   la place de $T_p|f$. ]

\mk
  
Si  $f$  est de degr\'e  $\leqslant k$  (comme polyn\^ome en $\D$), 
alors il en est de m\^eme de  $T_p|f$; on peut \'ecrire  $T_p|\D^k$ sous la forme
\begin{equation}\label{Tpfk1}
T_p|\D^k=\sum_{j=0}^k \mu_j \D^j,\quad {\rm avec } \; \mu_j\in \F_2.
\end{equation}
Supposons maintenant $k$ impair. Les formules (\ref{cnk}) et (\ref{ga2})
entra\^{\i}nent que
\begin{equation}\label{muj0}
j\not\equiv pk \pmod{8}\quad \Lr \quad \mu_j=0.
\end{equation}
 En particulier, on a $T_p(\cf_i) \subset \cf_j$ si $j \equiv pi \ $(mod $8)$.

L'opérateur de Hecke $T_p$ commute avec les opérations $f \mapsto
f^{2^s}$ de sorte que, si l'on conna\^it l'action de $T_p$ sur $\cf$,
par (\ref{fs}), on la conna\^it sur toutes les formes paraboliques.

\subsection{Nilpotence des op\'erateurs de Hecke modulo 2}
 
L'une des propri\'et\'es essentielles des op\'erateurs de Hecke modulo 2 
est qu'ils sont nilpotents (cf. par exemple \cite{Ha,Ono,Ser5}). Cela 
implique que, dans (\ref{Tpfk1}), le coefficient $\mu_k$ est nul. 
Par (\ref{Tpfk1}) et (\ref{muj0}), on a donc
pour tout $p$ premier $\geq 3$, et tout $k$ impair positif,
\begin{equation}\label{Tpfk}
T_p|\D^k=\sum_{\substack{j\equiv pk \hspace{-2mm}\pmod{8}\\1\leq j \leq k-2}} 
\mu_j \D^j ,\quad {\rm avec }\;\; \mu_j\in \F_2.
\end{equation}

\n {\it Exemples}:

\smallskip

{\rm (i)} $T_p|\D=0$ \ pour tout $p$ premier $ > 2$.

{\rm (ii)} Si $p\equiv 3 \pmod{8}$, on a $T_p|\D^3=\D;$ sinon, 
 $T_p|\D^3=0$. 

{\rm (iii)} Si $p\equiv 5 \pmod{8}$, on a $T_p|\D^5=\D;$ sinon, 
$T_p|\D^5=0$. 

{\rm (iv)} On a: 

$
\;\;    T_p|\D^7= 
\left\{
\begin{array}{llll}
   0  \quad {\rm si} \quad   p \equiv 1 \pmod{8}  & \! {\rm ou \ si } \quad  p \equiv -1 \pmod{16}\\ 
    \D^5 \  {\rm si}  \quad  p \equiv 3 \pmod{8}\\
   \D^3 \ {\rm si} \quad p \equiv 5 \pmod{8}\\
   \D \ \ {\rm si} \quad  p \equiv 7 \pmod{16}.
 \end{array}
\right. $

\subsection{L'ordre de nilpotence}\label{parindnil}
 
Par d\'efinition, {\it l'ordre de nilpotence} 
d'une forme modulaire $f\in \F_2[\D]$ est le plus petit entier 
$g=g(f)$ tel que, pour toute suite de  $g$  nombres premiers impairs  
$p_1,p_2,\ldots,p_g$, on ait 
$T_{p_1}T_{p_2}\ldots T_{p_g}|f=0$. 
[Comme les $T_p$ commutent entre eux, l'ordre dans lequel on \'ecrit 
les $T_{p_i}$ n'a pas d'importance. Noter aussi que l'on ne suppose 
pas que les  $p_i$ soient distincts.]
Lorsque $f$ = 0, on convient que $g(f)=-\iy$.

\mk

Nous d\'esignerons par $g(k)=g(\D^k)$ l'ordre de nilpotence de $\D^k$.  
Comme chaque $T_p$ abaisse le degr\'e en $\D$ d'au moins 2 unit\'es, on a
$g(k) \leq \frac{k+1}{2}\cdot$

\mk

Soit $p$ un nombre premier impair; il r\'esulte de la d\'efinition de 
l'ordre de nilpotence d'une forme modulaire $f\in\cf$ que l'on a 
\begin{equation}\label{gTpf}
g(f) \geq g(T_p|f) +1.
\end{equation}

\n {\it Exemples}:
 \begin{equation}\label{Tp3}
g(0)=-\infty,\ \ g(\D)=1,\ \ g(\D^3)= g(\D^3+\D)=2,\; \end{equation}
\begin{equation}\label{Tp5}g(\D^5)=g(\D^5+\D)=g(\D^5+\D^3)=g(\D^5+\D^3+\D)=2. \end{equation}

\section{Calcul des $T_p|\D^k :$ une r\'ecurrence lin\'eaire}\label{parThPol}

Soit $p$ un nombre premier $> 2$.
\begin{theoreme}\label{thHecPol}
Il existe un unique polyn\^ome  sym\'etrique $F_p(X,Y) \in \F_2[X,Y]$,
\begin{equation}\label{FpYX}
F_p(X,Y)=Y^{p+1}+s_1(X) Y^p +\ldots +s_p(X) Y +s_{p+1}(X)
\end{equation}
de degr\'e $p+1$ tel que
\begin{equation}\label{recTpDk}
T_p(\D^k) = \sum_{r=1}^{p+1} s_r (\D) \ T_p (\D^{k-r})
\end{equation}
pour tout $k \geq p+1$.
De plus, pour $1 \leq r \leq p+1$, $s_r(X)$ est une somme de mon\^omes en
$X$ dont les degr\'es sont congrus  \`a  $pr$ modulo $8$ et sont $\leq r$.
\end{theoreme}

\smallskip

\ni {\bf Esquisse de d\'emonstration}. On d\'efinit les  
$s_r(\D), \ 1 \leqslant i \leqslant p+1$, comme les fonctions
sym\'etriques \'el\'ementaires des $p+1$ s\'eries
$$f_0 = \D(q^p), \quad f_i = \D(z^iq^{1/p}), \quad i = 1,...,p,$$
o\`u $z$ est une racine primitive $p$-i\`eme de l'unit\'e dans une 
extension finie de $\F_2$. On d\'eduit (\ref{recTpDk}) de la formule:
$ T_p|\D^k = \sum_{i=0}^p \ (f_i)^k , \ \ k= 0,1,...$

\smallskip

\n {\it Exemples}
\footnote{Une table des polyn\^omes $F_p$ pour $p \leq 257$,
  calcul\'ee avec SAGE par Marc Del\'eglise, se trouve sur le site 
\tt{http://math.univ-lyon1.fr/$\sim$nicolas/polHecke.html}   }

\n Pour $p=3$ on a
\begin{equation}\label{F3XD}
F_3(X,Y)=(X+Y)^4 + XY =X^4+XY+Y^4.
\end{equation}
\n Vu (\ref{recTpDk}), cela donne un proc\'ed\'e de calcul 
des $T_3|\D^k$; si $t$ est une ind\'etermin\'ee, on a: 
$$
\sum_{k=1}^{\infty} T_3(\D^k)t^k = \frac{ \D t^3}{1+\D^3t+\D^4t^4}\cdot
$$
\n De m\^eme, pour $p=5$, on a:
\begin{equation}\label{F5XD}
F_5(X,Y)=(X+Y)^6 + XY = X^6 + X^4 Y^2 + X^2 Y^4 +X Y+Y^{6}
\end{equation}
et 
$$
\sum_{k=1}^{\infty} T_5(\D^k)t^k = \frac{ \D t^5}{1+\D^2t^2+\D^4t^4+\D^5t^5+\D^6t^6}\cdot
$$

\section{Les op\'erateurs de Hecke  $T_3$ et $T_5$} \label{parP3P5}

\subsection{Les nombres  $n_3(k), n_5(k)$ 
et $h(k)$}\label{parn3n5h}

Soit $k$ un nombre entier $\geq 0$. Ecrivons-le sous forme dyadique:
$\ds k=\sum_{i=0}^\iy \beta_i 2^i$ avec $\beta_i=0$ ou $1$. 
Posons:
$$ n_3(k)=\sum_{i=0}^\iy \beta_{2i+1} 2^i=
\sum_{\substack{i=1\\ i\;\rm{ impair}\;}}^\iy \beta_i 2^{\frac{i-1}{2}},\quad
n_5(k)=\sum_{i=0}^\iy \beta_{2i+2} 2^i=
\sum_{\substack{i=1\\ i\;\rm{ pair}\;}}^\iy \beta_i 2^{\frac{i-2}{2}},\quad
h(k)=n_3(k)+n_5(k).    
$$       

\n L'entier $h(k) $ est du m\^eme ordre de grandeur que $k^{1/2}$: si $k$ est
impair $> 0$ on a
$$ 
\frac 12 k^{1/2} \ < \ h(k)+1 \ < \  \frac 32 k^{1/2}.
$$ 

\n Notons que l'on a pour $\ell \geq 0$
$$
n_3(2\ell+1)=n_3(2\ell),\quad n_5(2\ell+1)=n_5(2\ell),\quad 
h(2\ell+1)=h(2\ell).
$$
Nous appellerons $[n_3(k),n_5(k)]$ le {\it code} du nombre $k$.
L'application   $k \mapsto [n_3(k),n_5(k)]$ est une bijection de l'ensemble des nombres 
impairs (resp. pairs) $\geq 0$ sur $\N^2$.

\subsection{Relation de domination}\label{pardominance}

Nous utiliserons la relation d'ordre suivante sur l'ensemble des nombres 
entiers naturels pairs (ou impairs): 
\begin{definition}\label{cord}
Si $k$ et $\ell$ ont m\^eme parit\'e, 
on dit que $\ell$ domine $k$ et on \'ecrit
$k\clet \ell$ ou $\ell \cget k$ 
si l'on a
$h(k) < h(\ell)$ ou bien $h(k)=h(\ell)$ et $n_5(k)< n_5(\ell)$.
La relation 
$k \cleq \ell$ d\'efinie par $k\clet \ell$ ou $k=\ell$,
est une relation d'ordre total sur l'ensemble des entiers pairs $($resp. impairs$) \geqslant 0$.
\end{definition}

\`A partir de maintenant, nous \'ecrirons une forme modulaire $f\in
\cf$, $f\neq 0$ sous la forme
\begin{equation}\label{fmi}
f=\D^{m_1}+\D^{m_2} \ldots +\D^{m_r} \;\;  {\rm avec }\;\;
m_1 \cget m_2 \cget \ldots \cget m_r.
\end{equation}

\subsection{La fonction $h$ pour les formes modulaires mod $2$}\label{parhformod2}

\begin{definition}\label{defhf}
Soit  $f\in\cf$. 

Si $f\neq 0$, on \'ecrit $f$ sous la forme (\ref{fmi}).
On dit que $m_1$ est {\rm l'exposant dominant} de $f$
et l'on d\'efinit $h(f)$ par  
$$h(f)=h(m_1)=\max_{1 \leq i \leq r} h(m_i). $$

Si $f=0$, on pose $h(f)=-\iy$.
\end{definition}

\subsection{Le cas de $T_3|f$} \label{parT3f}

\begin{proposition}\label{propP32}
Soit $f\in\cf$, $f\neq 0$  et soit $m_1$ son exposant dominant.

(i) On a $h(T_3|f) \leq h(f)-1=h(m_1)-1$.

(ii) Lorsque $n_3(m_1)\geq 1$, on a $h(T_3|f)=h(m_1)-1$ et 
l'exposant dominant 
de $T_3|f$ a pour code $[n_3(m_1)-1, n_5(m_1)]$.
\end{proposition}

\mk
\ni
{\bf D\'emonstration :} 
On consid\`ere d'abord le cas o\`u $f=\D^k$. 
On raisonne alors par r\'ecurrence sur $k$ en utilisant les relations 
(\ref{FpYX}), (\ref{recTpDk}) et (\ref{F3XD}).
La d\'emonstration est assez longue et technique.

\subsection{Le cas de  $T_5|f$} \label{parT5f}

\begin{proposition}\label{propP52}
Soit $f\in\cf$, $f\neq 0$  et soit $m_1$ son exposant dominant.

(i) On a $h(T_5|f) \leq h(f)-1=h(m_1)-1$.

(ii) Lorsque $n_5(m_1)\geq 1$, on a $h(T_5|f)=h(m_1)-1$ et 
l'exposant dominant 
de $T_5|f$ a pour code $[n_3(m_1), n_5(m_1)-1]$.
\end{proposition}

\mk
\ni
{\bf D\'emonstration :} 
M\^eme  m\'ethode que pour la proposition \ref{propP32};
on utilise (\ref{F5XD}) au lieu de (\ref{F3XD}).

\section{D\'etermination de l'ordre de nilpotence}
\label{demth1}

\begin{theoreme}\label{thmg=h}
Soit  $f\in\cf$, $f\neq 0$, que l'on \'ecrit
comme en {\rm (\ref{fmi})}.

(i) On a
\begin{equation}\label{T3T5fD}
T_3^{n_3(m_1)} T_5^{n_5(m_1)} |f=\D.
\end{equation}

(ii) La valeur de l'ordre de nilpotence $g(f)$ 
{\rm (cf. \S~\!\ref{parindnil})} est donn\'ee par
\begin{equation}\label{gf=hf}
g(f) = h(f)+1.
\end{equation}
\end{theoreme}

\mk
\ni
{\bf D\'emonstration :} 
(i) Soit $m$ l'exposant dominant de $\f=T_3^{n_3(m_1)} T_5^{n_5(m_1)}|f$. En appliquant $n_3(m_1)$ fois la proposition \ref{propP32} (ii) et
$n_5(m_1)$ fois la proposition \ref{propP52}~(ii), on voit que 
 $m$ a pour code $[0,0]$; comme $m$ est impair, on a $m=1$, d'o\`u $\f=\D$, ce qui 
d\'emontre (\ref{T3T5fD}). Notons que (\ref{T3T5fD}) implique
\begin{equation}\label{T3T5gf}
g(f) \geq n_3(m_1)+n_5(m_1)+1=h(m_1)+1=h(f)+1.
\end{equation}

(ii) Soit $d=\max(m_1,m_2,\ldots,m_r)$ le degr\'e de $f$; on va d\'emontrer 
(\ref{gf=hf}) par r\'ecurrence sur le nombre impair $d$.

Si $d=1,3$ ou $5$, (\ref{gf=hf}) r\'esulte de (\ref{Tp3}) et (\ref{Tp5}).

Soit $d \geq 7$ et supposons (\ref{gf=hf}) vraie pour toute 
forme de degr\'e $\leq d-2$. 
Pour $d \geq 7$, on a $h(d) \geq 2$ et la d\'efinition de l'exposant
dominant entra\^ine $h(f)=h(m_1) \geq h(d) \geq 2$. Par (\ref{T3T5gf}), on a
$g(f)\geq h(f)+1 \geq 3$; donc il existe des nombres premiers impairs 
$p_1,p_2,\ldots, p_s$ avec $s=g(f)-1 \geq 2$ et 
\begin{equation}\label{ming2}
T_{p_1} T_{p_2}\ldots T_{p_s}|f \neq 0.
\end{equation}

\n Posons $\f=T_{p_s}|f$, et calculons $g(\f)$. De (\ref{ming2}), on d\'eduit

\smallskip

$T_{p_1} T_{p_2}\ldots T_{p_{s-1}}|\f =T_{p_1} T_{p_2}\ldots T_{p_s}|f \neq 0$, 

\smallskip
\n ce qui implique $g(\f) \geq s$. Mais (\ref{gTpf}) entra\^ine
$g(\f)=g(T_p|f)\leq g(f)-1=s$. On en d\'eduit
\begin{equation}\label{ming4}
g(\f)= s =g(f)-1\geq 2.
\end{equation} 
Observons que (\ref{ming2}) et $s\geq 2$ entra\^inent $\f\neq 0$.
Par (\ref{Tpfk}), le degr\'e de $\f$ est $\leq d-2$;
on peut donc appliquer \`a  $\f$ l'hypoth\`ese de r\'ecurrence, ce 
qui donne $g(\f)=h(\f)+1$. En d\'esignant par $j$ l'exposant 
dominant de $\f$, 
avec (\ref{ming4}), il vient
\begin{equation}\label{ming5}
g(\f)=h(\f)+1=h(j)+1=s\geq 2.
\end{equation}
Soit $[u,v]$ le code de $j$,  avec $u\geq 0$, $v\geq 0$ et $u+v=s-1$.
En appliquant (i) \`a  $\f$ et en posant 
$q_1=q_2=\ldots =q_u=3$ et
$q_{u+1}=q_{u+2}=\ldots = q_{u+v}=5$, il vient
$$T_{q_1} T_{q_2}\ldots T_{q_{s-1}}|\f = 
T_{q_1} T_{q_2}\ldots T_{q_{s-1}}T_{p_s}|f = \D.$$
Posons $\psi=T_{q_{s-1}}|f$; on a
$$T_{q_1} T_{q_2}\ldots T_{q_{s-2}}T_{p_s}|\psi = 
T_{q_1} T_{q_2}\ldots T_{q_{s-1}}T_{p_s}|f = \D.$$ 
Cette formule montre que $g(\psi)\geq s$. 
Mais (\ref{gTpf}) entra\^ine
$g(\psi)=g(T_{q_{s-1}}|f)\leq g(f)-1=s$ \quad et \quad $g(\psi)=s$.

\smallskip 

\n Par (\ref{Tpfk}), le degr\'e de $\psi$ est  $\leq d-2$
et l'hypoth\`ese de r\'ecurrence 
donne $g(\psi)=h(\psi)+1$. On a ainsi
\begin{equation}\label{ming6}
g(\psi)= s =g(f)-1=h(\psi)+1.
\end{equation}
Par la proposition \ref{propP32}  (i) lorsque $q_{s-1}=3$, et par la
proposition \ref{propP52} (i) lorsque $q_{s-1}=5$, on a
$h(T_{q_{s-1}}|f) \leq h(f)-1$, d'o\`u, par (\ref{ming6}),
$$s-1=g(f)-2 = h(\psi)=h(T_{q_{s-1}}|f) \leq h(f)-1$$
ce qui implique $g(f) \leq h(f)+1$; vu  (\ref{T3T5gf}),
cela entra\^ine (\ref{gf=hf}).

\mk

\begin{corollaire}\label{corothm2}
Soit $f\in\cf$, $f\neq 0$, et soit $p$ un nombre premier tel que 
$p\equiv \pm 1 \pmod{8}$. Alors, on a 
\begin{equation}\label{corothm21}
g(T_p|f)\leq g(f)-2.
\end{equation}
\end{corollaire}

\mk
\ni
{\bf D\'emonstration :} 
On observe que, pour $p\equiv \pm 1 \pmod{8}$, on a $h(T_p|f) \equiv
h(f)\pmod{2}$,
ce qui, par le th\'eor\`eme \ref{thmg=h}, entraîne 
$g(T_p|f)\equiv g(f) \pmod{2}$.

\mk

\begin{corollaire}\label{corothm2bis}
Soit $f\in\cf$, $f\neq 0$. Si \ $T_3|f = T_5|f = 0,$ alors $f=\D$.
\end{corollaire}

\mk
\ni
{\bf D\'emonstration :} 
En effet, d'apr\`es (i), on a $n_3(m_1)=n_5(m_1)=0$, 
d'o\`u $m_1=1$ et $f=\D$.




\end{document}